
\input amstex
\documentstyle{amsppt}

\hyphenation{Hoch-schild}
                                        %
                                        %
\define\Bg{\goth g}                                          %
\define\Bgb{{\goth g}^{\bullet} }                                  %
\define\Lam{\wedge}                                        %
\define\Gam{\Gamma(M,\Lam^{\bullet}(TM))}                     %
\define\Gao{\Gamma(M,\Lam^{\bullet + 1}(TM))}                 %
\define\Cin{C^{\infty}(M)}                                   %
\define\Cbu{C^{\bullet}(A,A) }                               %
\define\Cbo{C^{\bullet + 1}(A,A)}                           %
\define\Ccu{C_{\bullet}(A,A)}                              %
\define\Ab{\Cal{A}^{\bullet}}                              %
\define\Bin{B_{\infty}}
\define\Gin{G_{\infty}}                     %
\define\Bvi{BV_{\infty}}                                     %
\define\Ain{A_{\infty}}                                       %
\define\Cinf{C_{\infty}}                                    %
\define\Din{D_{\infty}}                                 %
\define\Lin{L_{\infty}}                                     %
\define\La{\Lam ^{\bullet}(\text {Lie} (\Ab [1]^*))}        %


\topmatter
\title
Noncommutative differential calculus, homotopy BV algebras \\
and formality conjectures
\endtitle
\rightheadtext{Noncommutative differential calculus
and formality}
\author
D. Tamarkin and
B. L. Tsygan
\endauthor
\affil
Harvard University \\
Pennsylvania State University
\endaffil


\date
Received: 27.09.1999
\enddate
\thanks
The second author was partially supported by NSF Grant.
\endthanks
\dedicatory
In memory of Y. L. Daletski
\enddedicatory
\abstract
We define a notion of a strong homotopy BV algebra
and apply it to deformation theory problems.
Formality conjectures for Hochschild cochains
are formulated.
We prove several results supporting these conjectures.
\endabstract
\endtopmatter

\document
\centerline{Contents}
\vskip .2in
\roster
\item"0."  Introduction.
\item"1."  Homotopy Gerstenhaber algebras.
\item"2."  Homotopy BV algebras.
\item"3."  Formality theorems for Hochschild cochains.
\item"4."  Formality theorems for Hochschild chains.
\item"5."  Deformation complexes of bialgebras and homotopy BV operators.
\endroster


\head
0.\quad Introduction
\endhead


In non commutative geometry, a manifold $M$ is replaced by a possibly
non commutative algebra $A$ over a commutative ground ring $k$
(\cite{C}, \cite{G}). Many constructions of differential
calculus on $M$ (vector fields, polyvector fields, forms, etc.) have
their non commutative generalizations. For example, in case of vector
fields, for any algebra $A$ one can define the Lie algebra $Der (A)$
of derivations of $A$. If $A = C^{\infty}(M)$, one recovers the Lie
algebra of vector fields.

If one passes to the graded space $\Gamma (M, \Lam ^{\bullet}(TM))$
of polyvector fields on $M$, the non commutative analogue of this
construction is not just a graded space but a complex of Hochschild
cochains $C^{\bullet}(A,A)$ equipped with the Hochschild differential
$\delta :  C^{\bullet}(A,A) \rightarrow C^{\bullet + 1}(A,A)$ defined
for any algebra $A$ (cf.  \cite{CE}, \cite{Ger}). By definition,
$C^n(A,A)$ is the space of linear maps from $A^{\otimes n} $ to $A$.
If $A = C^{\infty}(M)$ then the subcomplex of those cochains which
are multidifferential operators has $\Gam$ as its cohomology
(\cite{HKR}). In this case, for us $\Cbu$ will always stand for this
subcomplex.

Now assume that a classical object has an algebraic structure (for
example, a structure of a graded algebra of some kind). Then one
could expect that its non commutative analogue (which is, as we see,
a complex) should have a structure of a differential graded algebra
of the same kind. We will see that this is true in many cases, and is
expected to be true in many more, if one replaces differential graded
algebras by more general objects, strong homotopy algebras (as
suggested by \cite{GJ1}, \cite{HJ}, \cite{K}; cf. also \cite{DT}).
The theory of strong homotopy structures was developed by Stasheff,
Moore and others, mainly for topological applications (cf., for
example, \cite{S}, \cite{LS}, \cite{Mo}). One can generalize the
notion of a differential graded associative algebra to define a
strong homotopy associative algebra (or $\Ain$ algebra); in the same
spirit one defines strong homotopy commutative algebras ($\Cinf$
algebras) and strong homotopy Lie algebras ($\Lin$ algebras). Below
we will have more examples of strong homotopy algebras.

A differential graded algebra is a strong homotopy algebra. However,
in a strong homotopy category morphisms are defined differently even
for ordinary algebras. A crucial property of strong homotopy
structures is that for them being quasi-isomorphic is an equivalence
relation. (One says that $\Cal{A}^{\bullet}_1$ is quasi-isomorphic to
$\Cal{A}^{\bullet}_2$ if there is a morphism from
$\Cal{A}^{\bullet}_1$ to $\Cal{A}^{\bullet}_2$ which induces an
isomorphism on cohomology).

The space $\Gao$ possesses a structure of a graded Lie algebra given
by the Schouten - Nijenhuis bracket $[\quad, \quad]_S$. On the non
commutative side, the Gerstenhaber bracket $[\quad, \quad]_G$ makes
$C^{\bullet + 1}(A,A)$ a differential graded Lie algebra for any
algebra $A$. By the formality theorem of Kontsevich, if $A = \Cin$,
then the differential graded Lie algebras
($\Cbo, \delta, [\quad, \quad]_G$) and \newline
($\Gao, \delta = 0, [\quad, \quad]_S $) are quasi-isomorphic as
strong homotopy Lie algebras, or $L_{\infty}$ algebras (cf.
\cite{K}).

        From this example one sees what should be the principles of
non commutative calculus:

\roster
\item"A)." An object from classical calculus should have its non
commutative analogue.
\item"B)." If there is an algebraic structure on a classical object,
the corresponding non commutative object should posess a similar
structure up to strong homotopy.
\item"C)." Thus, if $A=\Cin$, one gets two strong homotopy
structures, one coming from classical calculus and another from non
commutative one. Those two structures should be equivalent
(formality).
\endroster

As an algebraic structure gets richer, B), obviously, gets more
difficult to prove. On the other hand, C) becomes easier. Indeed, as
a general rule, the structure coming from non commutative calculus is
equivalent to a deformation of the classical structure; but the
richer the structure, the harder it is to deform. Let us illustrate
this principle by the following example from \cite{T}.

A graded space $\Cal{A}^{\bullet}$ is a Gerstenhaber algebra if two
operations are defined:
$$\cdot : \Cal{A}^{i} \otimes \Cal{A}^{j} \rightarrow \Cal{A}^{i + j}$$
and
$$[\quad,\quad] : \Cal{A}^{i+1} \otimes \Cal{A}^{j+1} \rightarrow
\Cal{A}^{i + j +1}$$
such that $\Cal{A}^{\bullet}$ is a graded commutative algebra with
respect to the product $\cdot$, $\Cal{A}^{\bullet + 1}$ is a graded
Lie algebra with respect to the bracket $[\quad,\quad]$, and the two
operations satisfy the Leibnitz identity
$$[a,b\cdot c] = [a,b] \cdot c + (-1)^{(|a|-1)|b|}b \cdot [a,c]$$
for any homogeneous elements $a,\; b, \;  c$ of $\Cal{A}^{\bullet}$.
It is well known that the Hochschild cohomology $HH^{\bullet}(A,A)$,
which is by definition the cohomology of the Hochschild complex
$\Cbu$, has a canonical structure of a Gerstenhaber algebra
(\cite{Ger}).

In \cite{T}, the following is proved.

\proclaim{\underbar{Theorem 0.1}}
For any algebra $A$ the space $\Cbu$ is a strong homotopy
Gerstenhaber algebra. Its underlying $\Lin$ structure is given by the
Gerstenhaber bracket. The induced structure of a Gerstenhaber algebra
on cohomology is given by the wedge product and the Schouten bracket.
\endproclaim

Cf. Section 1 for the definition of a strong homotopy Gerstenhaber
algebra. We call such algebras $G_{\infty}$ algebras.

\proclaim{\underbar{Theorem 0.2}}
For any algebra $A$,  the $G_{\infty}$ algebra $\Cbu$ is
quasi-isomorphic to a deformation of the $G_{\infty}$ algebra $HH^{\bullet}(A,A)$.
\endproclaim

\proclaim{\underbar{Theorem 0.3}}
If $M= \Bbb{R}^n$ and $A = \Cin$, the $G_{\infty}$ algebra $\Gam$ has
no non-trivial $G_{\infty}$ deformations. For any $M$, $\Cbu$ and
$\Gam$ are quasi-isomorphic as $G_{\infty}$ algebras.
\endproclaim
Theorem 0.3 generalizes the formality theorem of Kontsevich. We would
like to apply the above scheme to other constructions of classical
differential calculus.

We sketch the proofs of the above theorems in Section 3. For Theorem
0.1 a new version of the proof is presented, more streamlined as
compared to \cite{T}. This version is based on using Etingof-Kazhdan
dequantization instead of quantization, an idea suggested to us by P.
Etingof.

Now let us discuss non commutative analogues of other constructions
from calculus. A non commutative analogue of the space of
differential forms is well known to be the Hochschild chain complex
$C_{\bullet}(A,A)$ with the differential
$b: C_{\bullet}(A,A) \rightarrow C_{\bullet - 1}(A,A)$ (\cite{CE},
\cite{HKR}). One sees in this example that the principle B) above has
its obvious limitations: for example, the exterior product of forms
does not have an analogue on the Hochschild chain complex, unless $A$
is commutative. (A non commutative analogue of the product can be
restored if one passes from the algebra of forms to the algebra of
all differential operators on forms, cf. \cite{NT}). In this paper we
will mainly discuss not an algebra structure but a module structure
on forms.

There are two standard pairings between polyvector fields and forms:
for $\pi \in \Gamma(M, \Lam ^k (TM))$, one defines
$$i_{\pi} : \Omega ^{\bullet}(M) \rightarrow  \Omega ^{\bullet - k}(M)$$
and
$$L_{\pi} : \Omega ^{\bullet}(M) \rightarrow  \Omega ^{\bullet - k +1}(M)$$
where $i_{\pi}$ is the contraction of a form by a polyvector and
$$ L_{\pi} = [d,i_{\pi}]$$
The following identities hold (we denote the Schouten bracket
$[\quad,\quad]_S$ simply by $[\quad,\quad]$):
$$[L_{\pi}, L_{\varphi}] = L_{[\pi, \varphi]}; \quad
[i_{\pi}, L_{\varphi}] = i_{[\pi, \varphi]}; \quad
[i_{\pi}, i_{\varphi}] = 0
\tag0.1$$
as well as
$$i_{\pi \cdot \varphi} = i_{\pi} \cdot i_{\varphi};\quad
L_{\pi \cdot \varphi} =
L_{\pi} \cdot i_{\varphi} + (-1)^{|\pi|} i_{\pi} \cdot L_{\varphi}
\tag0.2$$
How to rephrase these identities in terms of Gerstenhaber algebras? 

A
Gerstenhaber module over a Gerstenhaber algebra $\Cal{A}^{\bullet}$
is a graded space $\Cal{M}^{\bullet}$ together with a structure of a
Gerstenhaber algebra on $\Cal{A}^{\bullet} \oplus \Cal{M}^{\bullet}$
such that $\Cal{A}^{\bullet}$ is a Gerstenhaber subalgebra and
$\Cal{M}^{\bullet}\cdot \Cal{M}^{\bullet} = [\Cal{M}^{\bullet} ,
\Cal{M}^{\bullet}] = 0$.

For any Gerstenhaber algebra $\Cal{A}^{\bullet}$ one constructs
canonically another Gerstenhaber algebra
$\Cal{A}^{\bullet}[\epsilon]$ as follows: let $\epsilon$ be a formal
parameter of degree one, $\epsilon ^2 = 0$. Define
$$\Cal{A}^{\bullet}[\epsilon]^n = \Cal{A}^{n} + \epsilon \Cal{A}^{n - 1}$$
Define the product $\star$ and the bracket $\{ \quad,\quad\}$ on
$\Cal{A}^{\bullet}[\epsilon]$ to be $k[\epsilon]$ - bilinear
operations such that
$$a\star b = ab + \epsilon (-1)^{|a|}[a,b]$$
and
$$\{a,b\} = [a,b]$$
for any homogeneous elements $a$ and $b$ of $\Cal{A}^{\bullet}$. We
will call the Gerstenhaber algebra $\Cal{A}^{\bullet}[\epsilon]$
 the canonical deformation of $\Cal{A}^{\bullet}$ (with the odd
parameter $\epsilon$). It is convenient to describe the formulas (0.1,0.2) above in the language of canonical deformations as follows:
\proclaim{\underbar{Lemma 0.4}}
Formulas
$$(\pi + \epsilon \varphi) \star \alpha = (-1)^{\pi}i_{\pi}\alpha$$
$$\{ \pi + \epsilon \varphi, \alpha \} =
L_{\pi}\alpha + i_{\varphi}\alpha$$
define a structure of a $\Gam [\epsilon]$ - module on
$\Omega ^{\bullet}(M)$.
\endproclaim

How to formulate a non commutative analogue of the above situation?
First of all, for any algebra $A$, for any $D \in C^k(A,A)$ one can
define operators
$$i_D : C_{\bullet}(A,A) \rightarrow   C_{\bullet - k}(A,A)$$
and
$$L_D : C_{\bullet}(A,A) \rightarrow   C_{\bullet - k + 1}(A,A);$$
operators $L_D$ define on $\Ccu$ a structure of a graded module over
the differential graded Lie algebra $\Cbo$. We show in subsection 1.2
how to extend the construction of the canonical deformation
$\Cal{A}^{\bullet}[\epsilon]$ to strong homotopy Gerstenhaber
algebras; more precisely, we will construct a canonical strong homotopy Gerstenhaber algebra structure on $\Cal{A}^{\bullet}[\epsilon]$ for any strong homotopy Gerstenhaber algebra $\Cal{A}^{\bullet}$.

\proclaim{\underbar{Conjecture 0.5}}
For any algebra $A$, on $\Cbu$ there is a canonical structure of a
$G_{\infty}$ module over the canonical deformation
$(\Cbu [\epsilon], \delta)$. This structure induces the structure of
a module over the differential graded Lie algebra $\Cbo$ given by
operators $L_D$. The induced structure of a Gerstenhaber module on
cohomology is as in Lemma 0.4.
\endproclaim

\proclaim{\underbar{Corollary 0.5.1}}
For any algebra $A$, on $\Cbu$ there is a canonical structure of a
$L_{\infty}$ module over the differential graded Lie algebra
$\Cbu [\epsilon]$ which extends the structure of a module over the
differential graded Lie algebra $\Cbo$ given by operators $L_D$.
\endproclaim

This last statement was proved in \cite{DT}.

One can prove a statement which is weaker than Conjecture 0.5:

\proclaim{\underbar{Theorem 0.6}}
For any algebra $A$, there is a $G_{\infty}$ algebra structure on
$\Cbu [\epsilon]$ such that $\Ccu$ is a $G_{\infty}$ module over
$\Cbu [\epsilon]$. The induced structure on $\Ccu$ of a module over
the differential graded Lie algebra $\Cbo$ is given by operators
$L_D$.  The induced structure of a Gerstenhaber module on cohomology
is as in Lemma 0.4.
\endproclaim

When $A = \Cin$, one can show that the $G_{\infty} $ algebra $\Cbu
[\epsilon]$ from Theorem 0.6 is quasi-isomorphic to the canonical
deformation $\Gam [\epsilon]$. Modulo $\epsilon$, at the level of
$\Lin$ algebras, this quasi-isomorphism is given by the formality
theorem of Kontsevich. Therefore one has two $G_{\infty}$ modules
over $\Gam [\epsilon]$: one is $\Omega ^{\bullet}(M)$ and the other
is $\Ccu$. One can prove as in \cite{T} that those $G_{\infty}$
modules are quasi-isomorphic.

By virtue of the formality theorem of Kontsevich, the differential
graded Lie algebra $\Gao$ is quasi-isomorphic to $\Cbo$ which acts on
$\Ccu$; therefore the latter is an $L_{\infty}$ module over  $\Gao$.

\proclaim{\underbar{Theorem 0.7}}
(Formality theorem for Hochschild chains). As $L_{\infty}$ modules
over $\Gao$, $\Ccu$ and $\Omega ^{\bullet}(M)$ are quasi-isomorphic.
\endproclaim

We show in \cite{Ts} that Theorem 0.7 allows to compute the
Hochschild homology of $A(\pi)$, where $\pi$ is a Poisson stricture
and $A(\pi)$ is the deformation of $\Cin$ given by the Kontsevich theorem
(\cite{K}). $A(\pi)$ is an algebra over ${\Bbb{C}}[[t]]$ whose
reduction modulo $t$ is isomorphic to $\Cin$. In particular, one gets

\proclaim
{\underbar {Corollary 0.7.1}} The space of traces on the deformed
algebra $A(\pi)$ is isomorphic to the space of
${\Bbb{C}}[[t]]$-valued distributions that are zero on all Poisson
brackets.
\endproclaim

Now we would like to extend the above constructions of non
commutative calculus to include the de Rham differential $d$. Recall
that a Batalin-Vilkovisky  (or BV) algebra is a Gerstenhaber algebra
$\Cal{A}^{\bullet}$ together with an operator $\Delta :
\Cal{A}^{\bullet} \rightarrow \Cal{A}^{\bullet - 1}$ such that
$\Delta ^2 = 0$, degree of $\Delta$ is $1$, and
$$\Delta (ab) = \Delta(a)b + (-1)^{|a|}a\Delta(b) + (-1)^{|a|} [a,b]$$
for any homogeneous elements $a$ and $b$ of $\Cal{A}^{\bullet}$.The
canonical deformation $\Ab [\epsilon]$ has a canonical BV operator
$\Delta = \partial / \partial \epsilon$.

One defines a BV module
$(\Cal{M}^{\bullet}, \Delta_{\Cal{M}^{\bullet}})$ over a BV algebra
$\Ab$ as a Gerstenhaber module $\Cal{M}^{\bullet}$ such that
$\Cal{A}^{\bullet} \oplus \Cal{M}^{\bullet}$ is a BV operator on the
corresponding Gerstenhaber algebra $\Ab \oplus \Cal{M}^{\bullet}$.

\proclaim
{\underbar{Lemma 0.8}} The de Rham complex
$(\Omega ^{\bullet}(M), d)$ is a BV module over the BV algebra
$(\Gam [\epsilon], \partial / \partial \epsilon)$.
\endproclaim

How to define a notion of a strong homotopy BV algebra (or
$BV_{\infty} $ algebra)? This question is of certain importance
because it is relevant for studying mirror symmetries, generalized
period maps for Calabi - Yau manifolds, Frobenius manifolds, and
related topics (\cite{BK}. \cite{B}, \cite{M}). In Section 2, we
propose an answer which is strongly suggested by our context. This
answer combines the standard technique of strong homotopy structures
with Kravchenko's construction from \cite{Ko}. A $\Bvi$ operator in
our sense is an odd differential operator of square zero on a certain
infinite dimensional graded space subject to additional properties.
One of those properties is a restriction on the order (and the
principal symbols) of its graded components, the other is that is has
to preserve some natural filtration.  A $BV_{\infty}$ algebra in our
sense is an algebra not over an operad but over a PROP.

Let us mention the reasons why we prefer our definition of a
$BV_{\infty}$ algebra. First, it is substantially more relaxed than
the standard definition using Koszul operads.  In the context of this
definition there is still a natural way to define a $\Bvi$ morphism,
such that being quasi-isomorphic is an equivalence relation. One can
prove an analogue of another standard property of strong homotopy
algebras, namely that any $\Bvi$ algebra is isomorphic to a direct
sum of a minimal $\Bvi$ algebra and a linearly contractible $\Bvi$
algebra (Section 2; compare to \cite{K}). Another reason is the
following. Note that the $G_{\infty} $ structure on $\Cbu$ which is
given by Theorem 0.1 is of a special kind: it arises from the chain
complex of a differential graded Lie bialgebra (cf. subsection 1.2).
For such algebras, the obstructions to being a $BV_{\infty}$ algebra
in our sense lie in the cohomology of the deformation complex of this
Lie bialgebra \cite{LR}.
Conjecturally, the constructions from
\cite{EK} extend to relate this cohomology to that of the deformation
complex of the Hopf algebra which is the canonical quantization of
the above Lie bialgebra; this last Hopf algebra is the cobar
construction of the algebra of Hochschild cochains equipped with the
product from \cite{GJ}, \cite{GV}; cf. also \cite{T}.

We show in Section 2 that for any $\Gin$ algebra $\Ab$, the canonical
deformation $\Ab [\epsilon]$ has a canonical $\Bvi$ structure.

\proclaim
{\underbar {Conjecture 0.9}} For any algebra $A$, the cyclic
differential $B$ on the Hochschild chain complex $\Ccu$ extends
canonically to a structure of a $\Bvi$ module over the $\Bvi$ algebra
$\Cbu [\epsilon]$, the canonical deformation of the Gerstenhaber
algebra $\Cbu$ from Theorem 0.1.
\endproclaim

If the above is true, then one easily proves the following

\proclaim
{\underbar {Corollary 0.9.1}} For any algebra $A$, the negative
cyclic complex \newline
$CC^-_{\bullet}(A) = (C_{\bullet}(A,A)[[u]], b+uB)$
has a canonical structure of an $\Lin$ module over the differential
graded Lie algebra
$(\Cbu [\epsilon][[u]], \delta + u\partial / \partial \epsilon)$.
Modulo $u$ this structure is the same as in Corollary 0.5.1.
\endproclaim

A statement which is very close to the above corollary was proven in
\cite{DT}. Note that the fact that was used in proving index theorems
(for example, in \cite{BNT}) is the existence of the fundamental
homomorphism $\chi$ from \cite{DT}, Section 8. This is a corollary of
0.9.1. Another, more effective way to construct this homomorphism is
contained in \cite{NT}, Section 4.

Now let us return to formality theorems for cyclic chains. Note that,
because of Theorem 0.3, for $A = \Cin$ $\Bvi$ algebras
$\Gam [\epsilon]$ and $\Cbu [\epsilon]$ are quasi-isomorphic.  Thus,
$\Ccu$ and $\Omega ^{\bullet} (M)$ are both $\Bvi$ modules over
$\Gam [\epsilon]$.

\proclaim
{\underbar {Conjecture 0.10}} For $A = \Cin$, $\Ccu$ and
$\Omega ^{\bullet} (M)$ are quasi-isomorphic as $\Bvi$ modules over
$\Gam [\epsilon]$.
\endproclaim

In particular, the following formality conjecture for cyclic chains
would be true (cf. \cite{Ts}):

\proclaim
{\underbar {Conjecture 0.11}} For $A = \Cin$, $\Ccu$ and
$\Omega ^{\bullet} (M)$ are quasi-isomorphic as $\Lin$ modules over
$\Gao [\epsilon]$.
\endproclaim

Let us mention a geometric application (cf. \cite{Ts}). Let $\pi$ be
a Poisson manifold $M$. By definition $\pi$ is regular if its
symplectic leaves form a foliation. Denote this foliation by
$\Cal{F}_{\pi}$. The tangent bundle to $\Cal{F}_{\pi}$ has $Sp (2n)$
as its structure group. Reducing the structure group to the maximal
compact subgroup, one can view $\Cal{F}_{\pi}$ as a $U(n)$-bundle;
let $\widehat{A}(\Cal{F}_{\pi})$ be the $\widehat{A} $ polynomial of
Chern classes of this bundle. The following is a corollary of the
conjecture above:

\proclaim
{\underbar {Conjecture 0.12}} The construction of
$\widehat{A}(\Cal{F}_{\pi})$ can be generalized to the non-regular
case.
\endproclaim

Now let us turn to another structure from classical calculus.
Assume that we are given an $n$-dimensional manifold $M$ with a
volume form $\Omega$. One defines the divergence operator
$$\Delta:  \Gamma(M, \Lam^{\bullet}(TM)) \rightarrow
\Gamma(M, \Lam^{\bullet - 1}(TM))$$
to be $\Delta = I^{-1} d I$ where
$I: \Gamma(M, \Lam^{\bullet}(TM)) \rightarrow \Omega^{n - \bullet} (M)$
sends a polyvector $\pi$ to the form $i_{\pi} \Omega$.
It is well known (cf., for example, \cite{BK}) that $\Delta$ is a BV
operator on the Gerstenhaber algebra $\Gam$.

As a non commutative analogue of the above situation, one can
consider an algebra $A$ with a trace $Tr$. We assume the pairing
$<a,b> = Tr(ab)$ to be non degenerate, meaning, essentially, that for
any Hochschild one-cochain there exists unique one-cochain $D^*$ such
that $<a,Db> = <D^*a, b>$ for all $a$ and $b$ in $A$. Strictly
speaking, we need a somewhat stronger condition: for any Hochschild
$m$-cochain $D$ there exists a unique $m$-cochain $D^*$ such that
$$Tr(a_0D(a_1, \dots, a_m) ) = Tr(a_m D^*(a_0, \dots, a_{m-1}))$$
for all $a_i \in A$ (cf. \cite{Sh}; ; note that our definition of
Hochschild cochains sometimes depends on an algebra $A$, like in the
case $A = \Cin$).

Given such a trace, one can define a map of complexes
$$J : \Cbu \rightarrow \Ccu ^* $$
and a differential
$$B_0: \Cbu \rightarrow C^{\bullet - 1}(A,A)$$
such that, if $B$ is Connes' cyclic differential on $\Ccu$, one has
$B^*J = JB_0$; cf. \cite{CFS}

\proclaim
{\underbar {Conjecture 0.13}} The operator $B_0$ extends to a $\Bvi$
structure on the $\Gin$ algebra $\Cbu$.
\endproclaim

        From this, along the lines of \cite{T}, one would be able to
obtain the cyclic formality conjecture of Kontsevich (cf. \cite{Sh}).

Conjecture 0.13 would imply a weaker form of Conjecture 0.9 which is
a cyclic analogue of Theorem 0.6. Indeed, one can apply Conjecture
0.13 to the algebra $R$ which is a trivial extension of $A$ by the
square zero ideal $A^*$ (the space dual to $A$, viewed as an
$A$-bimodule). The trace on $R$ is given by
$$Tr (a, \lambda) = \lambda (1)$$
for $a \in A$, $\lambda \in A^*$.

Finally, let us mention three topics of non commutative calculus that
are left out of this paper. First, one can try to prove another
version of Conjecture 0.9 using ``an infinitesimal analogue'' of the BV
operad. Second, a question very closely related to Conjecture 0.13 is
the study of algebraic structures on the Hochschild and cyclic
complexes of Hopf algebras. These complexes were defined recently by
Connes and Moscovici in \cite{CM}, \cite{CM1} with the motivation to
be able to prove index theorems by direct calculations (as was, in a different context, the
motivation for \cite{DT}, \cite{NT}). Third, as we mentioned
above, there is a non commutative analogue of the algebra
$D(\Omega ^{\bullet}(M))$ of differential operators on differential
forms (\cite{NT}). This is an $\Ain$ algebra whose cohomology is
$D(\Omega ^{\bullet}(M))$, with an $\Ain$ derivation acting on
cohomology as $[d_{DR}, \_]$ where $d_{DR}$ is the de Rham
differential. As a complex, this $\Ain$ algebra is equal to
$C_{\bullet}(\Cbu, \Cbu)$, the chain complex of the algebra of
cochains of $A$. The product on it is very closely related to the
product from \cite{GJ} and \cite{GV} that we are using in Section 3.
The first term of the $\Ain$ derivation is the cyclic differential
$B$ (the higher terms are very closely related to \cite{GJ1} and
\cite{HJ}).

The algebra $D(\Omega ^{\bullet}(M))$ has an additional structure, an
``order'' filtration $F_{\bullet}$ where, for a function $f$ and a
vector field $X$, operators $f$ and $i_X$ are of ``order'' zero and
operators $L_X$ and $df$ are of ``order'' one. Thus we have the
following algebraic structure: a graded algebra $D$ with an
increasing filtration $F_{\bullet}$ such that $gr_FD$ is commutative,
together with an extra derivation $d$ of degree one which satisfies
the transversality property that $dF_p$ is inside $F_{p+1}$ (as in the theory of variations of Hodge structures).

Such an algebraic structure admits a strong homotopy version (a
$D_{\infty}$ algebra). Conjecturally, for any $A$,
$C_{\bullet}(\Cbu, \Cbu)$ has a canonical $\Din$ structure. This
statement generalizes Theorem 0.1, as well as the main result of
\cite{NT}. It is closely related to the theorems and conjectures
above, as well as to constructions of \cite{B}. The exact nature of
these relationships is not yet clear.

{\bf Acknowledgements}. The authors are grateful to O.~Kravchenko, J. Stasheff , D.~Sullivan and A.~Voronov for many interesting discussions of homotopy BV algebras, and to 
J.~Stasheff for reading the preliminary version and making many helpful suggestions.
\subhead
1.\quad{\underbar{Homotopy Gerstenhaber algebras}}
\endsubhead
For a graded Lie algebra $\Bgb$ put
$$\Lam ^{k}\Bgb = S^k(\Bgb [-1]) $$
(as usual, for a complex,
$V^{\bullet}$, $V^{\bullet}[m]^i = V^{i + m}$). There is unique
bracket $[\quad,\quad]$ on $\bigoplus_{k\geq 0}\Lam ^{k}\Bgb$ such
that:

\roster
\item"1." The restriction of $[\quad,\quad]$ to
$\Lam ^{1}\Bgb = \Bgb$ is the Lie bracket on $\Bgb$, and
$[\Lam ^{0}\Bgb, \Lam ^{\bullet}\Bgb] = 0$.
\item"2."  With this bracket and with the exterior product,
$\Lam ^{\bullet}\Bgb$ becomes a Gerstenhaber algebra.
\endroster
If $\Bgb = \Bg$, the Lie algebra of a Lie group $G$, then this
Gerstenhaber algebra is the algebra of left invariant polyvector
fields on $G$.

We want to give a definition of a $\Gin$ structure on a graded space
$\Ab$. In what follows, we will take a liberty of talking about
algebras, not coalgebras; therefore, strictly speaking, our
definition will be correct only in the case when $\Ab$ is finite
dimensional. The definition we are about to give translates into a
collection of maps from $(\Ab )^*$ to $(\Ab) ^{* \otimes n}$ subject
to some identities ($*$ stands for dual); to make it correct, one has
to replace those maps by maps from $(\Ab) ^{ \otimes n}$ to $\Ab$,
subject to dual identities. Furthermore, all our strong homotopy
algebra structures will be described in terms of differential
operators on some algebra; that algebra will be always complete with
respect to a filtration, and we will always assume our operators to
be continuous in the corresponding topology.

For any graded space $V^{\bullet}$ consider the free graded Lie
algebra ${\text{Lie}}(V^{\bullet})$ generated by $V^{\bullet}$. Let
$\text{Lie}^k (V^{\bullet})$ be defined inductively by
$$
\text{Lie}^1 (V^{\bullet}) = V^{\bullet};\;
\text{Lie}^k (V^{\bullet}) = [V^{\bullet},
\text{Lie}^{k - 1}(V^{\bullet})]
$$
One has
$$ \text{Lie} (V^{\bullet}) =
\bigoplus_{k=1}^{\infty} \text{Lie}^k (V^{\bullet})$$
By ${\Cal{F}}^{\bullet}$ we will denote the filtration on
$ \bigoplus \Lam ^k (\text{Lie} (V^{\bullet}))$ induced by the
filtration ${\Cal{F}}^{\bullet} = \text{Lie}^{\geq k} (V^{\bullet})$,
and by $ \Lam ^{\bullet} (\text{Lie} (V^{\bullet}))$ the completion
of $\bigoplus \Lam ^k$ with respect to this filtration.

As usual, we define the grading on the dual space to a graded space
$V^{\bullet}$ by $(V^{\bullet})^{*i} = (V^{-i})^*$.

\definition
{\underbar {Definition 1.1}} A $\Gin$ algebra is a graded space $\Ab$
together with an operator $\partial$ of degree $1$ on
$\Lam (\text {Lie} ^{\bullet} (\Ab [1]^*))$ such that
$\partial ^2 = 0$ and $\partial $ is a derivation with respect to
both the product and the bracket.
\enddefinition

\remark
{\underbar {Remark \rom{1.2}}} Sometimes a different grading on $\La$ is
used, in which the bracket has degree $0$ and the product has degree
$1$. The reason for this is to make $\La$ an algebra over the operad
dual to $e_2$ where $e_2$ is the operad algebras over which are
Gerstenhaber algebras.
\endremark

A derivation $\partial $ from Definition 1.1 preserves two
multiplicative ideals of \newline
$\La$: $I_1$, generated by the
commutant of $\text {Lie}  (\Ab [1]^*))$, and \newline
$I_2 = \Lam ^{\geq 2} (\text {Lie} (\Ab [1]^*))$. The quotient of
$\La$ by $I_1$ is equal to \newline
$\Lam ^{\bullet}(\text{Lie}(\Ab [1]^*))$; since
$$ \Lam ^k (\Ab [1]^*) \simeq S^k (\Ab [1]^*[-1]) = S^k(\Ab [2]^*),$$
the differential induced by $\partial $ is by definition an $\Lin$
structure on the space $\Ab [1]$. On the other hand, the quotient of
$\La$ by $I_2$ is isomorphic to $\text {Lie}(\Ab [1]^*)$, thus
$\partial$ induces on it a differential which is by definition a
structure of a $\Cinf$
algebra on $\Ab$.

For a $\Gin$ algebra $\Ab$, a derivation $\partial$ of
$\Lam ^{\bullet}(\text{Lie}(\Ab [1]^*))$ is uniquely determined by its
restriction to $\Ab [1]^*$. The components of this restriction are
maps
$$\Ab [1] ^* \rightarrow
\text{Lie}^{k_1}(\Ab [1] ^*) \wedge \ldots \wedge
\text{Lie}^{k_n}(\Ab [1] ^*)$$
We view $ \text{Lie}(\Ab [1] ^*)$ as the subspace of primitive
elements of $U( \text{Lie}(\Ab [1] ^*) = T(\Ab [1]^*)$, therefore its
dual is a quotient of $T(\Ab [1])$. This quotient is taken by the
linear span of shuffle products of tensors of positive degree.

We denote by $m_{k_1, \ldots, k_n}$ the map conjugate to the map
above:
$$m_{k_1, \ldots, k_n}:
(\Ab)^{\otimes k_1} \otimes \ldots \otimes (\Ab)^{\otimes k_n}
\rightarrow \Ab$$
These operations have the following properties:
\roster
\item"1." Degree of $m_{k_1, \ldots, k_n}$ is equal to
$3-n-k_1 -\ldots -k_n$.
\item"2." Operations $m_{k_1, \ldots, k_n}(a_{11}, \ldots, a_{1,k_1};
\ldots ;a_{n1}, \ldots, a_{n,k_n})$ are Harrison \newline
cochains in every
group of variables $a_{i1}, \ldots, a_{i,k_i})$ (recall that a
Harrison cochain is a multi-linear map annihilating all shuffle
products).
\item"3." Operations $m_{k_1, \ldots, k_n}$ are invariant under
permutations of blocks \newline
$(a_{i1}, \ldots, a_{i,k_i})$ with appropriate
signs.
\endroster

The condition $\partial ^2 = 0$ now can be translated into a
quadratic identity for $m_{k_1, \ldots, k_n}$. In particular,
$\delta = m_1$ is a differential of degree $1$ on
$\Ab$; \newline
$a  b = (-1)^{|a|}m_2 (a,b)$ is a product of degree zero and,
more generally, operations $m_k$ define a $\Cinf$ structure on $\Ab$;
$[a,b] = (-1)^{|a| - 1}m_{1,1}(a,b)$ is a bracket of degree $-1$ and,
more generally, the operations $m_{1,\ldots,1}$ define a $\Lin$
structure on $\Ab [1]$; the operation $m_{1,2}(a;\;b,c)$ has degree
$-2$ and is up to a sign a homotopy for the Leibnitz identity, etc.

\example
{\underbar {Example 1.2}}
Let $\Ab$ be a Gerstenhaber algebra. Then one makes it a $\Gin$
algebra by defining the only non-zero operations
$m_{k_1, \ldots, k_n}$ to be
$m_{1,1}(a,b) = (-1)^{|a|}ab$, $m_2 (a,b) = (-1)^{|a|-1}[a,b]$.
\endexample

\example
{\underbar {Example 1.3}} Let $\Bgb$ be a graded Lie bialgebra
\cite{Dr}. On $\Lam^{\bullet}(\Bgb)$, let $\partial ^{\text {coLie}}$
be the derivation whose restriction to $\Bgb$ is the cobracket
$\Bgb \rightarrow \Lam ^2 (\Bgb)$. This is the standard cochain
differential of the Lie coalgebra $\Bgb$. If $\Bgb$ is free as a Lie
algebra, $\Bgb = \text{Lie}(V^{\bullet})$, then
$\partial ^{\text {coLie}}$ defines a $\Gin$ structure on the graded
space $\Ab = (V^{\bullet})^*[-1]$. More generally, if $\Bgb$ is a
differential graded Lie bialgebra, i.e. if a differential $\delta$ is
given on $\Bgb$ which is a derivation of the Lie algebra and a
coderivation of the Lie coalgebra, then we can extend $\delta$ to a
derivation of $\Lam^{\bullet} (\Bgb)$, and
$\partial = \delta + \partial ^{\text {coLie}}$ defines a $\Gin$
structure on $\Ab$.
\endexample

\subhead
1.1. \ A canonical deformation of a homotopy Gerstenhaber algebra
\endsubhead

First of all, note that for any graded Lie algebra $\Bgb$ the
Gerstenhaber algebra $\Lam ^{\bullet} (\Bgb)$ has a canonical BV
structure. Indeed, one can define $\Delta _{\text{Lie}}$ to be the
only BV operator satisfying
$$ \Delta _{\text{Lie}} (v \wedge w) = [v, w]$$
for $v,\;w$ in $\Bgb$. This operator is the usual chain differential
$\Lam ^{\bullet} (\Bgb) \rightarrow \Lam ^{\bullet - 1} (\Bgb)$.

Recall that for any Gerstenhaber algebra $\Ab$ one can define a
canonical deformation $(\Ab [\epsilon], \star, \{\quad,\quad\})$ (cf.
Introduction). Along with it one can define the trivial deformation
$\Ab [\epsilon]$ by extending the product and the bracket
$k[\epsilon]$ - bilinearly.

\proclaim
{\underbar {Lemma 1.2.1}} Let $\Delta$ be a BV operator on a
Gerstenhaber algebra $\Ab$. Then the operator $1 - \epsilon \Delta$
is an isomorphism between the trivial and the canonical deformations
of $\Ab$.
\endproclaim

The proof is straightforward.

Now, for a $\Gin$ algebra $\Ab$ let $\La [\epsilon]$ be the canonical
deformation of the Gerstenhaber algebra $\La $.
By Lemma above, the BV operator $\Delta _{\text{Lie}}$
defines an isomorphism  between the trivial and the canonical
$\Gin$ algebra structures on $\La [\epsilon]$. The derivation $\partial$ acts on the canonical
deformation, therefore it acts on the trivial deformation; this derivation
of the trivial deformation defines a $\Gin$ structure on
$\Ab [\epsilon]$; this is, by definition, the canonical deformation
of $\Ab$.

Let us give an explicit formula for the canonical deformation of a $\Gin$
structure, i.e. for the canonical $\Gin$ structure on $\Ab [\epsilon]$.

Let $D = [\partial, \Delta _{\text{Lie}}]$. This is a
derivation, and not just a differential operator of order two, precisely
because $\partial$ is a derivation with respect to the bracket. Then
the canonical deformation is defined by the derivation
$$\partial _{\text {can}} = \partial + \epsilon D
\tag1.1$$
of $\La [\epsilon]$. Strictly speaking, we need a derivation of
$\Lam ^{\bullet}(\Ab [\epsilon]^*)$, not of $\La [\epsilon]$. An easy
way to get one from another is to write down the
$k[\epsilon]$-multilinear operations $m_{k_1, \ldots,  k_n}$ defined
by $\partial _{\text {can}}$ and to notice that they satisfy all the
needed identities.

In the Example 1.3, note that any Lie bialgebra has a canonical
derivation $D : \Bgb \rightarrow \Lam ^2 (\Bgb) \rightarrow \Bgb$
which is the cobracket followed by the bracket. This is a derivation
and a coderivation. When extended to a derivation of
$\Lam ^{\bullet} (\Bgb)$, it coincides with the derivation $D$ above.

\proclaim
{\underbar {Lemma 1.2.3}} If $\Ab$ is a Gerstenhaber algebra (cf.
Example 1.2) then the above definition coincides with the one in
Introduction.
\endproclaim

\head
2. \quad Homotopy BV algebras
\endhead
 Recall that a BV operator $\Delta_{\text{Lie}}$ was defined in the
beginning of 1.2.

\definition
{\underbar {Definition 2.1}} A $\Bvi$ structure on a $\Gin$ algebra
$\Ab$ is an operator $\Delta$ on \newline
$\La$ such that:
\roster
\item "1." \hfill
$\Delta ^2 = 0; \;
\Delta = \partial + \Delta_1 + \Delta_3 + \ldots$\hfill\newline
where $\Delta _ i$ is an operator of degree $-i$ and $\partial$ is
the derivation defining the $\Gin$ structure on $\Ab$;
$\Delta _i (1) = 0$.
\item"2." $\Delta _1 - \Delta_{\text{Lie}} $ is a differential
operator of order $1$.
\item"3." $\Delta _{2i-1} $ is a differential operator of order $i$
for $i >1$.
\item"4." All operators $\Delta _k$ preserve the filtration
${\Cal{F}}^{\bullet}$.
\endroster
\enddefinition

We use the term differential operators in the sense of Grothendieck's
inductive definition. We can extend this definition as follows:

\definition
{\underbar{Definition 2.2}} Let
$f: R^{\bullet} \rightarrow S^{\bullet}$ be a morphism of graded
commutative algebras. A homogeneous operator
$D: \; R^{\bullet} \rightarrow S^{\bullet}$ is a differential
operator of order $k$ with respect to $f$ if the operator
$x \mapsto D(ax) - (-1)^{|D||a|}f(a)D(x)$ is of order $k-1$ with
respect to $f$ for any homogeneous $a$ in $ R^{\bullet}$; the zero
operator is of order $-1$ with respect to $f$.
\enddefinition

\definition
{\underbar{Definition 2.3}} A morphism of two $\Bvi$ algebras
$(\Ab, \; \Delta)$ and $({\Cal{B}^{\bullet}}, \Delta ')$ is an
operator
$$ F: {\Lam ^{\bullet}(\text {Lie}  ({\Cal {B}}^{\bullet} [1]^*))}
\rightarrow {\Lam ^{\bullet}(\text {Lie} ({\Cal {A}}^{\bullet} [1]^*))}$$
such that:
\roster
\item"1." \hfill $\displaystyle F = \sum _{i=0} ^{\infty}F_{2i}$
\hfill\newline
where $F_{2i}$ is of degree $-2i$.
\item "2." $F_{0}$ is a morphism of Gerstenhaber algebras.
\item "3." For $i>0$ $F_{2i}$ is a differential operator of order $i$
with respect to $F_0$ and $F_{2i}(1) = 0$.
\item"4." All the operators $F_k $ preserve the filtration
${\Cal{F}}^{\bullet}$.
\endroster
\enddefinition

A $\Bvi$ algebra is called minimal if the corresponding $\Gin$
algebra is minimal, i.e. if the differential $\delta = m_1$ is equal
to zero. A $\Bvi$ algebra $(\Ab, \Delta)$ is linear contractible if
the differential $\partial$ is induced by a contractible differential
on $\Ab$ and $\Delta = \partial + \Delta_{\text{Lie}}$.

\proclaim
{\underbar {Theorem 2.4}} Any $\Bvi$ algebra is isomorphic to direct
sum of a minimal $\Bvi$ algebra and a linear contractible $\Bvi$
algebra.
\endproclaim

\proclaim
{\underbar {Corollary 2.5}} If a there is a quasi-isomorphism
$\Cal{A}^{\bullet} \rightarrow \Cal{B}^{\bullet}$ then there is a
quasi-isomorphism $\Cal{B}^{\bullet} \rightarrow \Cal{A}^{\bullet}$.
\endproclaim

The proofs follow the usual scheme, as in \cite{K}, 4.5.

\remark
{\underbar {Remark \rom{2.6}}} Note that a $\Bvi$ structure in sense of
Definition 2.1 is in fact a homotopy BV operator from \cite{Ko}
acting on the algebra $\La$. In \cite{Ko}, however, the order of
$\Delta_{2i-1}$ is supposed to be $i+1$. Thus we can say that our
$\Bvi$ operator is Kravchenko's operator on $\La$ where the principal
symbol of $\Delta _1$ is that of $\Delta_{\text{Lie}}$ and the
principal symbols of $\Delta _i$ are equal to zero for all odd
$i > 1$.
\endremark

\proclaim
{\underbar {Proposition 2.7}} The canonical deformation of a $\Gin$
algebra has a canonical $\Bvi$ structure given by
$\Delta_{\text{can}} = \partial _{\text{can}} + \Delta _{\text {Lie}}
+ \frac{\partial}{\partial \epsilon}$
\endproclaim

The proof is straightforward.

\remark
{\underbar {Remark \rom{2.8}}} One can compare our definition to the usual
operadic definition. BV algebras are algebras over an operad which is
defined by relations of degree $\leq 2$ (the BV operad). One can
easily define Koszul operads in this context, following \cite{GK} and
\cite{Pr}. It is essentially proven in \cite{Get} and \cite{GJ} that
the BV operad is Koszul. Now, the standard operadic definition of a
$\Bvi$ algebra is that of an algebra over the bar construction of the
Koszul dual. This definition is a partial case of ours for which all
the operators $\Delta _i$ have to be derivations of $\La$ (they will
have automatically  to be derivations of the bracket).
\endremark

\head
3. \quad Formality theorems for Hochschild cochains
\endhead
In this section we outline the proofs of Theorems 0.1 - 0.3 from
Introduction.

{\underbar {Sketch of the proof of Theorem 0.1}} In order to show
that the complex of Hochschild cochains $\Cbu$ is a $\Gin$ algebra,
one shows first that it is a $B_{\infty}$ algebra, i.e. its cobar
construction is a differential graded Hopf algebra (see below). Then,
by a theorem of Etingof -  Kazhdan, one passes to the dequantization
of this Hopf algebra which is a Lie bialgebra. This Lie bialgebra is
cofree as a Lie coalgebra, and to such a bialgebra one can associate
a $\Gin$ algebra by a construction dual to that of Example 1.3. This
gives the required $\Gin$ structure on $\Cbu$.

More precisely, for a differential graded associative algebra
${\Cal {C}}^{\bullet}$, let $T({\Cal {C}}^{\bullet}[1])$ be the cobar
construction of ${\Cal {C}}^{\bullet}$, i.e. the differential graded
coalgebra with the comultiplication
$$\Delta(D_1\otimes \ldots \otimes D_n) = \sum _{k=1}^{n-1} (D_1
\otimes \ldots D_k) \bigotimes  (D_{k+1} \otimes \ldots \otimes D_n )
$$
and the differential
$$d(D_1\otimes \ldots \otimes D_n) =
\sum_{k=1}^{n-1}(-1)^{\sum_{i\leq k}(|D_i|+1)-1} D_1 \otimes \ldots
D_k D_{k+1} \otimes \ldots \otimes D_n + $$
$$\sum _{k=1}^{n}(-1)^{\sum_{i < k}(|D_i|+1)} D_1 \otimes \ldots
\delta D_i \otimes \ldots \otimes D_n $$
Here $\delta$ is the differential in ${\Cal {C}}^{\bullet}$ and $D_i$
are homogeneous elements of ${\Cal {C}}^{\bullet}$.

\definition
{\underbar {Definition 3.1}} We say that ${\Cal {C}}^{\bullet}$ is a
$\Bin$ algebra if there is a multiplication on the cobar construction
of ${\Cal {C}}^{\bullet}$ which makes it a differential graded
bialgebra. We require this multiplication to preserve the filtration
$f_p = T^{\leq p}({\Cal {C}}^{\bullet}$).
\enddefinition

Consider the following example (\cite{GJ}, \cite{GV}). Let
${\Cal{C}}^{\bullet}$ be the Hochschild cochain complex of an algebra
$A$. The algebra $A$ may be itself graded. For cochains $D_i$ in
$C^{d_i}(A, A)$ define a cochain
$$
\multline
D_0\{D_1, \ldots , D_m\}(a_1, \ldots, a_n) =\\
=\sum (-1)^{ \sum_{q <p \leq m}(|a_{i_q}|-1)(|D_p| -1)}  D_0(a_1,
\ldots ,a_{i_1} ,\\ D_1 (a_{i_1 + 1}, \ldots ),\ldots ,
D_m (a_{i_m + 1}, \ldots ) , \ldots)
\endmultline
$$
(the total degree $|D_i|$ of $D_i$ is its degree of homogeneity plus
 $d_i$). Let $m$ be the two-cochain of $A$ defined by
$$m(a,b) = (-1)^{|a|}ab$$
Put
$$D\smile E = (-1)^{|D|}m\{D,\;E\}; $$
$$ [D,\;E] = D\{E\}-(-1)^{(|D|-1)(|E|-1)}E\{D\};\quad \delta D = [m,D]$$
This defines respectively the cup product, the Gerstenhaber bracket
and the Hochschild differential on $\Cbu$. The space $\Cbu$ with the
differential $\delta$ and the cup product is a differential graded
associative algebra, the space $\Cbo$ with the bracket and the
differential is a differential graded Lie algebra, and the above
operations induce on $HH^{\bullet}(A,A)$ a structure of a
Gerstenhaber algebra.

To define a multiplication on $T(\Cbu)$ note that, since it has to be
compatible with the comultiplication, it has to be determined by its
composition with the projection from $ T(\Cbu)$ to $\Cbu$. We define
this projection of the product
$(D_1 \otimes \ldots D_m) \bullet (E_1 \otimes \ldots E_n)$ to be
equal to $D\{E_1, \ldots, E_n\}$ if $m=1$ and $0$ if $m > 1$.

\proclaim
{\underbar {Lemma 3.2}} The product $\bullet$ defines on the
differential graded algebra $(\Cbu, \smile, \delta)$ a $\Bin$
structure.
\endproclaim
Proof: cf. \cite{GJ}, \cite{GV}.

For a $\Bin$ algebra ${\Cal{C}}^{\bullet}$, the bialgebra
$T({\Cal{C}}^{\bullet})$ has an antipode $S$ which can be defined
explicitly:
$$\multline
\quad \quad S(D_1\otimes \ldots \otimes D_n) + D_1\otimes \ldots
\otimes D_n = \\
\sum_{0<k_1<n} (D_1\otimes \ldots \otimes D_{k_1})
\bullet (D_{k_1 + 1} \otimes \ldots \otimes D_{n}) - \\
\sum_{0<k_1<k_2<n} (D_{1} \otimes \ldots \otimes D_{k_1}) \bullet
(D_{k_1 + 1} \otimes \ldots \otimes D_{k_2}) \bullet (D_{k_2 + 1}
\otimes \ldots \otimes D_{n}) + \ldots
\endmultline
$$
Thus $T({\Cal{C}}^{\bullet}[1])$ becomes a differential graded Hopf
algebra. By the result from section 2.4
of \cite{EK}, suitably adapted to our situation,
this Hopf algebra admits a canonical dequantization which
is a graded Lie bialgebra. Since the construction of a dequantization
is natural, the differential on $T({\Cal{C}}^{\bullet}[1])$ defines a
differential on this Lie bialgebra. One can show that this bialgebra
is cofree as a Lie coalgebra, with the space of cogenerators
${\Cal{C}}^{\bullet}[1]$. Applying the construction from Example
1.3 to the dual Lie bialgebra, we get the required $\Gin$ structure
on ${\Cal{C}}^{\bullet}$. Note that for all $\Gin$ algebras which are
obtained as in Example 1.3 all the operations $m_{k_1, \ldots , k_n}$
are zero if $n >2$. Therefore, in particular, the induced $\Lin$
structure on ${\Cal{C}}^{\bullet}[1]$ does not have non-zero higher
brackets.

This proves Theorem 0.1. Theorem 0.2 follows easily from Theorem 2.4,
or rather from its version for $\Gin$ algebras. To prove Theorem 0.3,
one notices that the deformation complex of any $\Gin$ algebra $\Ab$
is the complex of derivations of the Gerstenhaber algebra $\La$, with
the differential $[\partial, \_]$. For $M = {\Bbb{R}}^n$ and
$\Ab = HH^{\bullet}(A,A) = \Gam$, the cohomology of this complex
could be computed explicitly: there is a spectral sequence converging
to it, with the first term equal to the Poisson cohomology (similar to
the one defined in \cite{Br}) of the odd symplectic space
$T^*(M) = {\Bbb{R}}^{n|n}$. Since this cohomology is essentially the de
Rham cohomology of $M$, one concludes that in case of
$M = {\Bbb{R}}^n$ any deformation of the $\Gin$ algebra $\Gam$ must
be trivial. Thus, by virtue of Theorem 0.2, Theorem 0.3 is true for
$M = {\Bbb{R}}^n$.  For general $M$, one needs an additional argument
with Gelfand - Fuks cohomology as in \cite{K}, \cite{T}.

\head
4. Formality theorems for Hochschild chains
\endhead
{\underbar {Sketch of the proof of Theorem 0.6}} For an algebra $A$,
consider the dual space $A^*$ as an $A$ - bimodule. Let
$R = A \dot{+} A^*$ be the algebra with the product
$$ (a+\lambda)(b+\mu) = ab + a\mu + \lambda b $$
for $a,\;b \in A,\; \lambda,\; \mu \in A^*$.
Any Hochschild cochain is a multi-linear form on $ A \dot{+} A^*$; it
can be decomposed into components which are linear forms on
$A^{\otimes p} \otimes (A^*)^{\otimes q}$ with values in $A$ or in
$A^*$. The full Hochschild complex contains a subcomplex of those
cochains whose components are given, respectively, by linear maps
$A^{\otimes p} \rightarrow A^{\otimes q} \otimes A$ or
$A^{\otimes p} \otimes A \rightarrow A^{\otimes q}$. Our definition
of tensor powers of $A$ may depend on the ring $A$; for example, if
$A = \Cin$ then
$$ A^{\otimes p} = \text{jets}_{\text{diagonal}}C^{\infty}(M^p) $$
Define a grading on $R$ by $R^0 = A$ and $R^{-1} = A^*$.Then the
complex $C^{\bullet}(R,R)$ becomes graded; all its components of
degree $<-1$ are equal to zero and
$$C^{\bullet}(R,R)^{-1} = C_{\bullet}(A,A)^*$$
(the dual complex).
The subcomplex $C^{\bullet}(R,R)^{0}$ consists of cochains whose
components (see above) are all zero except for $p=q=1$ or $p=q=0$.
The cochains with $p=q=1$ form a subcomplex; denote it by
$C^{\bullet}(R,R)^{1}_1.$ The quotient of $C^{\bullet}(R,R)^{0}$ by
this subcomplex is exactly $\Cbu$. It is easy to see that the exact
sequence
$$0 \rightarrow C^{\bullet}(R,R)^{1}_1 \rightarrow
C^{\bullet}(R,R)^{0} \rightarrow \Cbu \rightarrow 0$$
splits and that $C^{\bullet}(R,R)^{1}_{1}$ is quasi-isomorphic to
$C^{\bullet - 1}(A,A)$. Thus, the complex $C^{\bullet}(R,R)^{0}$ is
quasi-isomorphic to $\Cbu [\epsilon]$.

Next step is to show that the operations $D_0\{D_1, \ldots , D_m\}$
from Section 3 all preserve the grading on $C^{\bullet}(R,R)$. After
that, applying the results of Section 3, one sees that
$C^{\bullet}(R,R)^0$ is a $\Gin$ algebra and $C^{\bullet}(R,R)^{-1}$
is a $\Gin$ module over it. Since every $\Gin$ module admits a dual,
we see that $\Ccu$ is a $\Gin$ module over $C^{\bullet}(R,R)^0$ and
therefore over some odd deformation $\Cbu [\epsilon]$ of the $\Gin$
algebra $\Cbu$.

This implies Theorem 0.6. After that, using the same reasoning as in
the end of Section 3, one deduces Theorem 0.7.

\head
5.  Deformation complexes of bialgebras and homotopy BV operators
\endhead
Let $\Bgb$ be a graded Lie bialgebra. Put (cf. \cite{LR})
$$C^{p-1,q-1}(\Bgb) = \text{Hom}(\Lam^p(\Bgb), \Lam^q(\Bgb))$$
for $p,\;q\geq 1$. Recall that, in our notation,
$\Lam^q(\Bgb)=S^p(\Bgb [-1])$; for \newline
$D \in C^{p-1,q-1}(\Bgb)$ put
$|D| = 2p-2+$ (homogeneity degree of $D$). For example, if
$\Bgb = \Bg ^0 = \Bg$ is a Lie bialgebra, then $C^{p-1,q-1}(\Bgb)$ is
concentrated in degree $p+q-2$. There are two differentials
$$d^{\text{Lie}}: C^{p-1,q-1}(\Bgb) \rightarrow C^{p,q-1}(\Bgb)$$
and
$$d_{\text{coLie}}: C^{p-1,q-1}(\Bgb) \rightarrow C^{p-1,q}(\Bgb)$$
The first is the standard cochain differential of the Lie algebra
$\Bgb$ with coefficients in the module $\Lam^{q}(\Bgb)$, the second
is a dual differential for the coalgebra $\Bgb$ (or the differential
$d^{\text{Lie}}$ for the dual Lie algebra $(\Bgb)^*$). If $\delta$ is
a differential of the bialgebra $\Bgb$ which is a derivation with
respect to both structures, then it induces a differential
$$\delta: C^{p-1,q-1}(\Bgb) \rightarrow C^{p-1,q-1}(\Bgb)$$
The above three differentials commute with each other. We put
$d=\delta+d^{\text{Lie}}+d_{\text{coLie}}$. The total complex
$(C^{\bullet}, d)$ is called the deformation complex of the
differential graded Lie bialgebra $\Bgb$.

Let us interpret this complex in terms of differential operators on
the algebra $\Lam^{\bullet}(\Bgb)$. For any such differential
operator of order $p$, we say that it is of pure order $p$ if its
restriction to $\Lam ^{<p}(\Bgb)$ is zero. Let ${\Cal{D}}^p$ be the
space of differential operators of pure order $p$. Define the Poisson
bracket (cf. \cite{KS})
$$\{ \quad,\quad \}: {\Cal{D}}^p \otimes {\Cal{D}}^q \rightarrow
{\Cal{D}}^{p+q-1}$$
to be the component of pure order $p+q-1$ of the commutator of
operators.

Clearly, an operator of pure order $p$ is uniquely determined by its
restriction to $\Lam ^p (\Bgb)$, so we identify ${\Cal{D}}^p$ with
$C^{p-1,\bullet}(\Bgb)$. Note that the derivations $\delta$ and
$\partial ^{\text{coLie}}$ from Example 1.3 are in ${\Cal{D}}^1$ and
the chain differential $\Delta _{\text{Lie}}$ is in ${\Cal{D}}^2$.

\proclaim
{\underbar {Lemma 5.1}} Under the identification of $C^{\bullet}
(\Bgb)$ with the space of differential operators, one has
$$\delta = \{\delta, \_\};\;
\partial^{\text{Lie}} = \{\Delta _{\text{Lie}}, \_\};\;
\partial_{\text{coLie}} = \{\partial^{\text{coLie}},\_\}$$
where the left hand sides refer to the three differentials on
$C^{\bullet} (\Bgb)$ and the right hand sides to differential
operators on $\Lam ^{\bullet}(\Bgb)$.
\endproclaim

The proof is straightforward.

If a Lie bialgebra $\Bgb$ has a filtration preserved by the bracket,
the cobracket and the differential, then by the filtered deformation
complex of $\Bgb$ we will mean the subcomplex of $C^{\bullet}(\Bgb)$
consisting of maps preserving the induced filtration on
$\Lam ^{\bullet}(\Bgb)$.

Recall from Section 3 that a $\Bvi$ operator on the cochain complex
$\Cbu$ is a differential operator on $\Lam (\Bgb)$
$$\Delta =
\delta + \partial^{\text{coLie}} + \Delta_1 + \Delta_3 + \ldots$$
satisfying
$$\Delta ^2 = 0,$$
where $\Bgb = \Bgb(A)$ is the Etingof - Kazhdan dequantization of the
Hopf algebra $T(\Cbu [1])$. The order of
$\Delta _1 - \Delta _{\text{Lie}}$ has to be equal to $1$, the order
of $\Delta_ {2i - 1}$ to $i$ for $i > 1$. All $\Delta _i$ have to
preserve the increasing filtration on $\Lam (\Bgb)$ which is induced
by a certain filtration on $\Bgb$. (The filtration is increasing, not
decreasing as in Section 3, because we pass to the dual bialgebra).

\proclaim
{\underbar {Theorem 5.2}} All obstructions to existence of a $\Bvi$
structure on $\Cbu$ lie in the cohomology of the filtered deformation
complex of the differential graded Lie bialgebra $\Bgb (A)$.
\endproclaim

{\underbar{Proof}} Rewrite $\Delta$ as
$$\Delta = \sum_{p \geq 0} \Delta (p),$$
where $\Delta (0) = \delta + \partial^{\text{coLie}}$ and
$\Delta (p)$ is the sum of components of pure order $i-p$ in
$\Delta _{2i-1}$ for $p>0$. Solving the equation $\Delta ^2 = 0$ is
equivalent to solving a chain of equations
$\{\Delta(0), \Delta (p)\} = c_p$ where $c_p$ is a known cochain.

\proclaim
{\underbar {Corollary 5.3}} The first obstruction to existence of a
$\Bvi$ operator on $\Cbu$ is the canonical derivation $D$ of
$\Bgb (A)$ (the cobracket followed by the bracket).
\endproclaim

For a graded associative bialgebra $H^{\bullet}$, there is a parallel
version of a deformation complex (cf. \cite{GS}). Put
$$C^{p-1, q-1}(H^{\bullet}) =
\text{Hom}((H^{\bullet})^{\otimes p}, (H^{\bullet})^{\otimes q})$$
For $D \in C^{p-1, q-1}$ set
$|D| = p+q-2+(\text{homogeneity degree of D})$. Define
$$\partial ^{\text {alg}}: C^{p-1, q-1}(H^{\bullet}) \rightarrow
C^{p, q-1}(H^{\bullet})$$
$$\partial _{\text {coalg}}:C^{p-1, q-1}(H^{\bullet}) \rightarrow
C^{p-1, q}(H^{\bullet})$$
to be respectively the Hochschild cochain differential of the algebra
$H^{\bullet}$ with coefficients in the bimodule
$(H^{\bullet})^{\otimes p}$ and the dual differential for the
coalgebra $H^{\bullet}$ (or, in other words, the differential
$\partial ^{\text {alg}}$ for the dual algebra $(H^{\bullet})^*$). If
$\delta$ is a differential on $H^{\bullet}$ which is a derivation and
a coderivation, then it induces a differential
$$\delta: C^{p-1, q-1}(H^{\bullet}) \rightarrow
C^{p-1, q-1}(H^{\bullet})$$
Put
$d = \delta + \partial ^{\text {alg}} + \partial _{\text {coalg}}$.
The total complex $C^{\bullet}(H^{\bullet})$ with the differential
$d$ is called the deformation complex of $H^{\bullet}$.  If
$H^{\bullet}$ is a filtered bialgebra then the filtered deformation
complex consists of those maps which preserve the induced filtration
on tensor powers.

\proclaim
{\underbar {Conjecture 5.4}} There is a canonical quasi-isomorphism
between $C^{\bullet}(\Bgb)$ and \newline
$C^{\bullet}(H^{\bullet})$ where
$\Bgb$ is a Lie bialgebra and $H^{\bullet}$ is its canonical
quantization of Etingof - Kazhdan.
\endproclaim

\proclaim
{\underbar {Conjecture 5.5}} The image of the canonical derivation
$D$ under the  canonical dequantization of Etingof - Kazhdan is
$\frac{1}{2} \text{log}(S^2)$ where $S$ is  the antipode.
\endproclaim

The deformation complex of a Lie bialgebra has the usual algebraic
structure, namely the Lie bracket (in our grading it has degree
zero). As we see, it also has an associative product coming from its
identification with the space of differential operators. What are the
analogues of these structures for Hopf algebras, and how to translate
a $\Bvi$ operator on $\Cbu$ into their terms? The product
$C^{p-1,q-1} \otimes C^{p'-1,q'-1} \rightarrow C^{p+p'-1,q+q'-1}$
exists; it is the analogue of the commutative multiplication of
degree two on the deformation complex of a bialgebra (the
multiplication of principal symbols).

It is more or less clear that a Hopf algebra analogue of a $\Bvi$
operator should be multiplicative, not a ``connection'' like $\Delta$
but rather its ``holonomy''.


\Refs
\widestnumber\key{GDT1}

\ref
  \key B
  \by S. Barannikov
  \paper Generalized periods and mirror symmetry in dimension $n>3$ \newline
     \paperinfo   AG/9903124
 \endref

\ref
  \key BK
  \by S. Barannikov and M. Kontsevich
  \paper Frobenius manifolds and formality of Lie algebras
        of polyvector fields
  \paperinfo   AG/9710032
 \endref

\ref
  \key BV
  \by I.A. Batalin and G.S. Vilkovisky
  \paper Gauge algebra and quantizations
  \jour Phys. Letters
  \vol 102 B
  \yr 1981
  \pages 27-31
\endref

\ref
  \key BNT
  \by P. Bressler, B. Nest and B. Tsygan
  \paper Riemann-Roch theorems via deformation quantization, I
  \paperinfo   AG/9904121
 \endref

\ref
  \key BNT1
  \by P. Bressler, B. Nest and B. Tsygan
  \paper Riemann-Roch theorems via deformation quantization
  \paperinfo   AG/9705014
 \endref

\ref
  \key Br
  \by J.L. Brylinski
  \paper Differential complex for Poisson manifolds
  \jour JDG
  \vol 28
  \yr 1988
  \pages 93-114
\endref

\ref
  \key CE
  \by A. Cartan, and S. Eilenberg
  \paper Homological Algebra
  \paperinfo Princeton Univ. Press (1956)
\endref

\ref
  \key C
  \by A. Connes
  \paper Noncommutative Differential Geometry
  \paperinfo Publ. Math., IHES, vol. 62
\endref

\ref
  \key CFS
  \by A. Connes, M. Flato and D. Sternheimer
  \paper Closed star products and cyclic homology
  \jour Letters in Math. Physics
\vol 24
\yr 92
\pages 1-12
\endref

\ref
  \key CM1
  \by A. Connes and H. Moscovici
  \paper Cyclic homology and Hopf algebras
  \paperinfo QA/9904154
\endref

\ref
  \key CM
  \by A. Connes and H. Moscovici
  \paper Cyclic homology, Hopf algebras and modular theory
  \paperinfo  \newline
QA/9905013
  \jour Letters in Math. Phys.
  \vol 48
  \issue 1, to appear
\endref

\ref
  \key DT
  \by Y. Daletski and B. Tsygan
  \paper Operations on Hochschild and cyclic complexes
  \jour this issue
  \vol
  \yr
\pages
\endref

\ref
  \key Dr
\by V. Drinfeld
\paper Hamiltonian structures on Lie groups, Lie bialgebras and the geometric meaning of the classical Yang-Baxter equation
\jour Soviet Math. Dokl
\vol 27
\yr 1983, nr. 1
 \endref

\ref
  \key EK
\by P. Etingof and D. Kazhdan
\paper Quantization of Lie bialgebras, II
\jour Selecta Mathematica, New series
\vol 4
\yr 199\pages 213-231
 \endref

\ref
  \key FT1
  \by B. Feigin and B. Tsygan
  \paper Additive $K$-theory
  \paperinfo Springer Lect. Notes, \underbar{1289} (1987), 66--220
\endref

\ref
  \key GDT1
  \by I. Gelfand, Y. Daletski and B. Tsygan
  \paper On a variant of noncommutative geometry
  \jour Soviet Mat. Dokl.
  \vol 40
  \yr 1990
  \pages 422--426
\endref

\ref
  \key GD
  \by I. Gelfand and I. Dorfman
  \paper  Hamiltonian operators and associated algebraic structures
\jour I. M. Gelfand's Collected Papers
  \paperinfo vol. 1 (1987), 707--721
\endref

\ref
  \key G
  \by I. Gelfand and M. Naimark
  \paper On the imbedding of normed rings
        into the rings of operators in Hilbert space
  \jour Mat. Sb.
  \vol 12
  \yr 1943
  \pages 197--213
\endref

\ref
  \key Ger
  \by M. Gerstenhaber
  \paper The cohomology structure of an associative ring
  \jour Ann. Math
  \vol 78
  \yr 1963
  \pages 59--103
\endref

\ref
  \key GS
  \by M. Gerstenhaber and S. Schack
  \paper Algebras, bialgebras, quantum groups, and algebraic deformations
  \jour Contemporary Mathematics
  \vol 134
  \yr 1992
  \pages 51--92
\endref

\ref
  \key GV
  \by M. Gerstenhaber and A. Voronov
  \paper Homotopy $G$-algebras and moduli space operad
  \jour Int. Math. Res. Notices
  \yr 1995
\issue 2
\pages 141--153
\endref

\ref
  \key Ge
  \by E. Getzler
  \paper Cartan homotopy formulas
        and the Gauss Manin connection in cyclic homology
  \jour Israel Math. Conf. Proc.
\vol 7
\yr 1993
\endref

\ref
  \key Get
  \by E. Getzler
  \paper  Batalin-Vilkovisky algebras
        and two-dimensional topological field theories
\jour CMP
\vol 159
\yr 1994
\pages 265--185
\endref

\ref
  \key GJ
  \by E. Getzler and J. Jones
  \paper Operads, homotopy algebras
        and iterated integrals for double loop spaces
  \paperinfo hep-th/9403055
\endref

\ref
  \key GJ1
  \by E. Getzler and J. Jones
  \paper $\Ain$ algebras and the cyclic bar complex
\jour Illinois J. of Math.
\vol 34
\yr 1990
\pages 256--283
\endref

\ref
  \key GJP
  \by E. Getzler, J. Jones and S. Petrack
  \paper Differential forms on loop spaces and cyclic bar complex
\jour Topology
\vol 30
\yr 1991
\pages 339--371
\endref

\ref
  \key GK
  \by V. Ginzburg and M. Kapranov
  \paper Koszul duality for operads
\jour Duke Math. J.
\vol 76
\yr 1994
\issue 1
\pages 203--273
\endref

\ref
  \key Go
  \by T. Goodwillie
  \paper Cyclic homology, derivations and the free loop space
  \jour Topology
  \vol 24
  \yr 1985
  \pages 187--215
\endref

\ref
  \key HKR
  \by G. Hochschild, B. Kostant, A. Rosenberg
  \paper Differential forms on regular affine algebras
  \jour Trans. AMS
  \vol 102
  \yr 1962
  \pages 383--408
\endref

\ref
  \key HJ
  \by C.E. Hood and J. Jones
  \paper Some algebraic properties of cyclic homology groups
\jour K-Theory
\vol 1
\yr 1987
\pages 361-384
\endref

\ref
  \key K
  \by M. Kontsevich
  \paper Deformation quantization of Poisson manifolds, I
  \paperinfo   QA/9709040
 \endref

\ref
  \key KS
  \by Y. Kosmann-Schwarzbach
  \paper From Poisson algebras to Gerstenhaber algebras
  \jour Ann. Inst. Fourier, Grenoble
  \vol 46
  \yr 1996
  \pages 1243-1274
\endref

\ref
  \key Ko
  \by O. Kravchenko
  \paper Deformations of Batalin - Vilkovisky algebras
  \paperinfo   QA/9903191
 \endref

\ref
  \key LS
  \by T. Lada and J.D. Stasheff
  \paper Introduction to sh Lie algebras for physicists
  \jour Int. J. Theor. Phys.
  \vol 32
  \yr 1993
  \pages 1087-1103
\endref

\ref
  \key LR
  \by P. Lecomte and C. Roger
  \paper Modules et d\'{e}formations des big\`{e}bres
  \jour CRAS
  \vol 310
  \yr 1990
  \pages 405-410
\endref

\ref
  \key M
\by Yu.I. Manin
\paper Three constructions of Frobenius manifolds: a comparative study \newline
\paperinfo QA/98010006
 \endref

\ref
  \key Mo
  \by J.C. Moore
  \paper Alg\`ebre homologique et homologie des espaces
  \inbook S\'{e}minaire Henri Cartan, 12\`{e}me ann\'{e}e, 1959-60
  \publ W.A. Benjamin
  \publaddr New York - Amsterdam
  \yr 1967
  \endref

\ref
  \key Mo1
  \by J.C. Moore
  \paper Differential homological algebra
  \inbook Actes du CIM (Nice, 1970)
  \pages 335-339
  \publ Gauthiers-Villard
  \publaddr Paris
  \yr 1970
\endref

\ref
  \key NT
  \by R. Nest. B. Tsygan
  \paper On the cohomology ring of an algebra
  \jour Advances in Geometry, Birkh\"{a}user
  \vol 172
  \yr 1998
  \pages 337-370
\endref

\ref
  \key Pr
  \by S.B. Priddy
  \paper Koszul resolutions
  \jour Trans. AMS
  \vol 152
  \yr 1970
  \pages 39-60
\endref

\ref
  \key R
  \by G. Rinehart
  \paper Differential forms on general commutative algebras \newline
  \jour Trans. AMS
  \vol 108
  \yr 1963
  \pages 195--222
\endref

\ref
  \key Sh
  \by B. Shoikhet
  \paper On the cyclic formality conjecture
  \paperinfo math.QA/9903183
 \endref

\ref
  \key S
  \by J. Stasheff
  \paper Homotopy associativity of H-spaces
  \jour Trans. Amer. Math. Soc.
  \vol 108
  \yr 1963
  \pages 293-312
\endref

\ref
  \key T
  \by D. Tamarkin
  \paper Another proof of Kontsevich formality conjecture
  \paperinfo   QA/9803025
 \endref

\ref
  \key Ts
  \by B. Tsygan
  \paper Formality conjectures for chains
  \paperinfo   QA/9904132
 \endref

\ref
  \key Z
  \by B. Zwiebach
  \paper Closed string field theory: quantum action
        and the B-V master equation
  \jour Nuclear Physics
  \vol B390
  \yr 1993
  \pages 33-152
\endref

\endRefs
\enddocument